\documentclass[12pt]{article}
\usepackage[T1]{fontenc}
\usepackage[ansinew]{inputenc}
\usepackage{amsmath}
\usepackage{graphicx}
\title{Rigidity for nearly umbilical hypersurfaces in space forms}
\author{Xu Cheng and Detang Zhou\thanks{Both authors are  partially supported by CNPq and Faperj  of Brazil.}}
\usepackage{amsmath, amssymb}
\usepackage{amsthm}
\makeatletter

\@addtoreset{equation}{section} \makeatother
\newtheorem{thm}{Theorem}[section]

\theoremstyle{remark}

\begin{document}
\date{}
\maketitle
 \begin{abstract}
In \cite{P}, Perez proved  some $L^2$ inequalities for closed convex hypersurfaces immersed in the Euclidean space $\mathbb{R}^{n+1}$, more generally, for closed hypersurfaces with non-negative Ricci curvature, immersed in an Einstein manifold.  In this paper, we discuss   the rigidity of these inequalities when the ambient manifold is $\mathbb{R}^{n+1}$,   the hyperbolic space $\mathbb{H}^{n+1}$,  or  the closed hemisphere $\mathbb{S}_+^{n+1}$. We also obtain a generalization of  the Perez's theorem to the hypersurfaces without the hypothesis of non-negative Ricci curvature. 
\end{abstract}

\baselineskip=14pt
\section {Introduction}
In this paper,  we suppose 
 $\Sigma$ is a smooth connected oriented closed (i.e. compact and without boundary) hypersurface immersed in an $(n+1)$-dimensional Riemannian manifold $(M, \widetilde{g})$ with induced metric $g$. Recall that
 $\Sigma$ is called totally umbilical if its second fundamental form  $A$  is multiple  of its metric  $g$ at every  point  of $ \Sigma$, that is,  $A=\frac{H}{n}g$ on $\Sigma$.  Here,  $A$ is  defined by   $A(X,Y)=-\left<\widetilde{\nabla}_XY,\nu\right>$, where $\nu$ denotes the outward unit normal to $\Sigma$, $X, Y\in T\Sigma$, $ \widetilde{\nabla}$ denote the Levi-Civita connection of $(M, \widetilde{g})$.  $H=\text{tr}A$ denotes the mean curvature of $\Sigma$, which is the trace of $A$.  A classical theorem in  differential geometry states that a closed  totally umbilical surface in the Euclidean space $\mathbb{R}^3$ must be a round sphere $\mathbb{S}^2$  and its second fundamental form $A$ is a constant multiple of its metric $g$. This theorem also holds for higher dimensional cases.
There are  various generalizations  of this theorem (for instance, cf.  a survey \cite{R}). 

 In 2005,  De Lellis and M\"uller \cite{dLM} considered a stability of the above theorem and   proved that if $\Sigma\subset \mathbb{R}^3$ is a closed connected surface with normalized  area $4\pi$,  then
$$ ||A - \text{Id}||_{L^2(\Sigma)}\leq  C ||A-\frac{\mathrm{tr}A}{2}\text{Id}||_{L^2(\Sigma)} ,$$
 where  $C$ is a universal constant.

 Recently, D. Perez \cite{P} generalized  the inequality of De Lellis and M\"uller to convex hypersurfaces. He proved that 
 \begin{thm} (\cite{P})\label{thm-1}
Let $\Sigma$ be a smooth, closed and connected hypersurface in $\mathbb{R}^{n+1}, n\geq 2,$ with induced Riemannian $g$ and non-negative Ricci curvature, then 
\begin{equation}\label{ine-i2}\int_{\Sigma}|A-\frac1n\overline{H}g|^2\leq \frac{n}{n-1}\int_{\Sigma}|A-\frac{H}{n}g|^2,
\end{equation}
and equivalently 
\begin{equation}\label{ine-i-02}\int_{\Sigma}(H-\overline{H})^2\leq \frac{n}{n-1}\int_{\Sigma}|A-\frac{H}{n}g|^2,
\end{equation}
where $\overline{H}=\frac1{\text{Vol}_n(\Sigma)}\int_{\Sigma}H$. In particular, the above estimate holds for smooth, closed hypersurfaces which are the boundary of a convex set in $\mathbb{R}^{n+1}$.

\end{thm}
For a closed hypersurface in the Euclidean space $\mathbb{R}^{n+1}$,  it is known that $Ric\geq 0$ is equivalent  to $A\geq 0$ (that is, $\Sigma$ is convex) (see its proof, for instance,  in  \cite{P},  Pg.48).  So this implies that Theorem \ref{thm-1} holds for closed convex hypersurfaces.
In \cite{P}, the author also showed that the constants in  inequalities (\ref{ine-i2}) and (\ref{ine-i-02})  are sharp, and that  without assumption on nonnegativity of the Ricci curvature,  the inequality cannot hold with a universal constant. Moreover,  as pointed out by De Lellis and Topping \cite{dLT},  Perez's theorem holds even for the closed hypersufaces with non-negative Ricci curvature when the ambient space is Einstein.  Indeed this can be seen by an observation on the Codazzi's equation satisfied by hypersurfaces in an Einstein manifold (more precisely, see equation (\ref{eq-ip-1}) of this paper). Hence a slight modification of the  proof of Perez \cite{P} gives

 \begin{thm}  \label{thm-2}Let $(M^{n+1}, \widetilde{g}), n\geq 2,$ be an Einstein manifold.
Let $\Sigma$ be a smooth, closed and connected hypersurface immersed in $M$ with  non-negative Ricci curvature, then 
\begin{equation}\label{ine-i3}\int_{\Sigma}|A-\frac1n\overline{H}g|^2\leq \frac{n}{n-1}\int_{\Sigma}|A-\frac{H}{n}g|^2,
\end{equation}
and equivalently 
\begin{equation}\label{ine-i-03}\int_{\Sigma}(H-\overline{H})^2\leq \frac{n}{n-1}\int_{\Sigma}|A-\frac{H}{n}g|^2,
\end{equation}
where $\overline{H}=\frac1{\text{Vol}_n(\Sigma)}\int_{\Sigma}H$. 

\end{thm}

In this paper, we will study the rigidity of Perez's inequalities (Theorem \ref{thm-1} and Theorem \ref{thm-2}).  We study what happens if  (\ref{ine-i2}) and (\ref{ine-i-02}), or (\ref{ine-i3}) and (\ref {ine-i-03}) hold as an equality.   When the ambient spaces are $\mathbb{R}^{n+1}$, $\mathbb{H}^{n+1}$ and the closed hemisphere $\mathbb{S}_{+}^{n+1}$, we  prove that

\begin{thm}\label{thm-3} Assume   $M^{n+1}(c), c=0,-1,1,$ are  the Euclidean space $\mathbb{R}^{n+1}$,  the hyperbolic space $\mathbb{H}^{n+1}(-1)$, and the closed hemisphere $\mathbb{S}_{+}^{n+1}(1)$, respectively.
Let $\Sigma$ be a smooth, closed and connected hypersurface immersed in $M$ with  non-negative Ricci curvature. Then 
\begin{equation}\label{eq-i-1} \int_{\Sigma}|A-\frac1n\overline{H}g|^2=\frac{n}{n-1}\int_{\Sigma}|A-\frac{H}{n}g|^2,
\end{equation}
and  equivalently
\begin{equation}\label{eq-i-01}\int_{\Sigma}(H-\overline{H})^2=\frac{n}{n-1}\int_{\Sigma}|A-\frac{H}{n}g|^2,
\end{equation}
 hold if and only if $\Sigma$ is a totally umbilical hypersurface, where $\overline{H}=\frac1{\text{Vol}_n(\Sigma)}\int_{\Sigma}H$, that is, $\Sigma$ is a distance sphere $S^n$ in $M^{n+1}( c )$. 
\end{thm}

Here a distance sphere $S^n$ in a complete Riemannian manifold $M^{n+1}$ is defined as  the set of points in $M$ which have the same distance from a fixed point in $(M, \widetilde{g})$. It is known that a closed totally umbilical hypersurface in a space form is a distance sphere  (especially, a  distance sphere in $\mathbb{R}^{n+1}$ is a round sphere) and its  second fundamental form $A$ is a constant multiple of its metric.   Theorem \ref{thm-3} discuss  rigidity of the above known fact.

If the ambient manifold $M$ is any Einstein manifold, we cannot expect the equalities like (\ref{eq-i-1})  and (\ref{eq-i-01}) hold for a closed totally umbilical hypersurface in $M$. For instance, the complex projection space $\mathbb{CP}^{\frac{n+1}{2}}$ has no closed totally umbilical hypersurfaces. On the other hand,  the constants in Theorem \ref{thm-2}  are  also  sharp when the ambient Einstein manifold is a space form besides $\mathbb{R}^{n+1}$,  which was proved by A. Ju\'arez \cite{J}.

We further study the general case in which  the hypersurfaces has no assumption on its Ricci curvature. In this case, although an inequality with universal constant cannot hold as demonstrated in \cite{P},   we may still obtain    quantitative $L^2$ inequalities and discuss the rigidity  of these inequalities. Precisely, we prove that
 \begin{thm} \label{thm-4} Let $(M^{n+1}, \widetilde{g}), n\geq 2,$ be an Einstein manifold.
Let $\Sigma$ be a smooth, connected and closed  hypersurface immersed in $M$ with induced metric $g$.  Then 
\begin{equation}\label{ine-i-4}\int_{\Sigma}|A-\frac{\overline{H}}{n} g|^2\leq \frac{n}{n-1}(1+\frac{K}{\eta_1})\int_{\Sigma}|A-\frac{H}{n}g|^2,
\end{equation}
and equivalently
\begin{equation}\label{ine-i-04}\int_{\Sigma}(H-\overline{H})^2\leq\frac{n}{n-1}(1+\frac{nK}{\eta_1})\int_{\Sigma}|A-\frac{H}{n}g|^2,
\end{equation}
where $\eta_1$ is the first nonzero eigenvalue of the Laplacian operator on $\Sigma$, $K\geq 0$ is a  nonnegative constant such that the Ricci curvature of $\Sigma$  satisfies $\text{Ric} \geq -K$.

Moreover,  when  $M^{n+1}$ is  the Euclidean space $\mathbb{R}^{n+1}$, the hyperbolic space $\mathbb{H}^{n+1}(-1)$, or the closed hemisphere $\mathbb{S}_{+}^{n+1}(1)$, the equality in (\ref{ine-i-4}) (or equivalently (\ref{ine-i-04})) holds if and only if $\Sigma$ is  a totally umbilical hypersurface,  that is, $\Sigma$ is a distance sphere $S^n$ in $M^{n+1}$.
\end{thm}

Note that take $K=0$ in Theorem \ref{thm-4}, we obtain Theorem \ref{thm-2} and Theorem \ref{thm-3}. Hence Theorem \ref{thm-4} generalizes Theorem \ref{thm-2} and Theorem \ref{thm-3}.
We take  the method used in \cite{P} to prove inequalities (\ref{ine-i-4}) and   (\ref{ine-i-04}) in Theorem \ref{thm-4}. On proof of   Theorem \ref{thm-3} and the case of equalities in Theorem \ref{thm-4}, 
the  methods used in \cite{dLT} and \cite{GWX} could not be applied.  However we can use the submanifold theory to prove them. 

It is worth  mentioning  that there is a similar phenomenon in Riemannian geometry.  The classical Schur lemma states that the scalar curvature of an Einstein manifold of dimension $n \geq 3$ must be constant.  In \cite{dLT},  De Lellis and Topping  discussed the stability and rigidity of this assertion and proved an almost-Schur lemma.  Recently, the first author of the present paper \cite{C} generalized the almost-Schur lemma to the case in which the Ricci curvature has no non-negative lower bound and  in \cite{C2} obtained a  generalization of almost-Schur lemma for symmetric $(2,0)$-tensors and made some applications.

\section{Proof of theorems}

Throughout the paper,   we use the same notation to denote a symmetric $(2,0)$-tensor and its corresponding $(1,1)$-tensor.  For example, $\text{Ric}$  denotes the symmetric   $(2,0)$ Ricci tensor of $\Sigma$ and its corresponding $(1,1)$ Ricci tensor.
\bigskip

{\it Proof of Theorem \ref{thm-4}}. Take a local orthonormal frame $\{e_i\}, 1\leq i\leq  n,$ on $\Sigma$, which can be extended to a local orthonormal frame $\{e_i\}, 1\leq i\leq n+1$ on $M$,  where $e_{n+1}=\nu$.  $h_{ij}=A(e_i,e_j), 1\leq i,j\leq n$.
The Codazzi equations are
$$h_{ij,k}-h_{ik,j}=\widetilde{R}_{n+1 ijk}, 1\leq i, j, k\leq n$$
where $\widetilde{R}$ denotes the Riemannian curvature tensor on $(M,\widetilde{g})$.

Take $j=i$ in the Codazzi equations and take the sum of the index $i$ from $1$ to $n$. Then
\begin{equation}
H_k=\displaystyle\sum_{i=1}^nh_{ii,k}=\displaystyle\sum_{i=1}^nh_{ik,i}+\displaystyle\sum_{i=1}^n\widetilde{R}_{n+1 iik}=\displaystyle\sum_{i=1}^nh_{ik,i}+\widetilde{R}_{n+1 k}.\nonumber
\end{equation}
Since $M$ is Einstein, $\widetilde{R}_{n+1 k}=\frac{R}{n}\widetilde{g}_{n+1 k}=0.$ We have
\begin{align}\label{eq-ip-1}\nabla H=\textrm{div} A,
\end{align}
where $\textrm{div} A=\displaystyle\sum_{i=1}^n(\nabla_{e_i}A) (e_i,e_k)=\displaystyle\sum_{i=1}^nh_{ik,i}e_k$.   

Denote by $\text{\r A}=A-\frac{H}{n}g$  the traceless tensor of $A$.
 Then
$$\textrm{div}  \text{\r A}=\textrm{div}  A-\textrm{div} (\frac{H}{n}g)=\textrm{div}  A-\frac{\nabla H}{n}.$$
So by (\ref{eq-ip-1}),   
\begin{equation}\label{eq-p-1}\nabla H=\frac{n}{n-1}\textrm{div}  \text{\r A}.
 \end{equation} 
Let $\phi$ be the unique solution of the Poisson equation on $\Sigma$:
\begin{equation}\label{eq-p-2}
\Delta \phi=H-\overline{H}, \quad \int_{\Sigma}\phi=0.
\end{equation}
If $\phi\equiv \text{Constant}$ on $\Sigma$, inequality (\ref{ine-i-04}) obviously holds. So now we suppose $\phi$ is not identically zero on $\Sigma$. 
We have
\begin{align}\label{ine-p-1}
\int_{\Sigma}|H-\overline{H}|^2&=\int_{\Sigma}(H-\overline{H})\Delta \phi=-\int_{\Sigma}\left<\nabla H,\nabla\phi\right>\nonumber\\
&=-\frac{n}{n-1}\int_{\Sigma}\left<\text{div} \text{\r A},\nabla\phi\right>\nonumber\\
&=\frac{n}{n-1}\int_{\Sigma}\left< \text{\r A},\nabla^2\phi\right>\nonumber\\
&=\frac{n}{n-1}\int_{\Sigma}\left< \text{\r A},\nabla^2\phi-\frac1n(\Delta\phi)g\right>\nonumber\\
&\leq\frac{n}{n-1}\left(\int_{\Sigma}| \text{\r A}|^2\right)^{\frac12}\left[\int_{\Sigma}|\nabla^2\phi-\frac1n(\Delta\phi)g|^2\right]^{\frac12}\nonumber\\
&=\frac{n}{n-1}\left(\int_{\Sigma}| \text{\r A}|^2\right)^{\frac12}\left[\int_{\Sigma}|\nabla^2\phi|^2-\frac1n\int_{\Sigma}(\Delta\phi)^2\right]^{\frac12}
\end{align}
Applying the Bochner formula to $\phi$ ,
integrating, and applying the Stokes'  formula, we have
\begin{equation}\label{eq-p-3}\int_{\Sigma}|\nabla^2\phi|^2=\int_{\Sigma}(\Delta\phi)^2-\int_{\Sigma}\text{Ric}(\nabla\phi,\nabla\phi).
\end{equation}
By (\ref{ine-p-1}) and (\ref{eq-p-3}), we have 
\begin{equation}
\int_{\Sigma}|H-\overline{H}|^2\leq\frac{n}{n-1}\left(\int_{\Sigma}| \text{\r A}|^2\right)^{\frac12}\left[\frac{n-1}{n}\int_{\Sigma}(\Delta\phi)^2-\int_{\Sigma}\text{Ric}(\nabla\phi,\nabla\phi)\right]^{\frac12}
\end{equation}
Since $\text{Ric}\geq -(n-1)K, K\geq 0$, 
\begin{equation}\int_{\Sigma}\text{Ric}(\nabla \phi, \nabla \phi)\geq -(n-1)K\int_{\Sigma}|\nabla \phi|^2.\nonumber
\end{equation}
\begin{align}\label{ine-pp-1}
\int_{\Sigma}|H-\overline{H}|^2&\leq\frac{n}{n-1}\left(\int_{\Sigma}| \text{\r A}|^2\right)^{\frac12}\left[\frac{n-1}{n}\int_{\Sigma}(\Delta\phi)^2+(n-1)K|\nabla \phi|^2\right]^{\frac12}\nonumber\\
&=\sqrt{\frac{n}{n-1}}\left(\int_{\Sigma}| \text{\r A}|^2\right)^{\frac12}\left[\int_{\Sigma}(\Delta\phi)^2+nK|\nabla \phi|^2\right]^{\frac12}
\end{align}
Let $\eta_1$ denote the first nonzero eigenvalue of the Laplace operator  on $\Sigma$, i.e., 
\[\eta_1=\displaystyle\inf\{ \frac{\int_M |\nabla \varphi|^2}{\int_M \varphi^2};  \varphi\in C^{\infty}(M) \textrm{ is not identically zero and }\int_M \varphi=0 \}.
\]   
 We have
\begin{align}
\int_{\Sigma} |\nabla \phi|^2 
&=-\int_{\Sigma}\phi\Delta \phi=-\int_{\Sigma}\phi(H-\overline{H})\nonumber\\
&\leq \left(\int_{\Sigma}\phi^2\right)^{\frac12}\left[\int_{\Sigma}(H-\overline{H})^2\right]^{\frac12}\nonumber\\
&\leq \left(\frac{\int |\nabla \phi|^2}{\eta_1}\right)^{\frac12}\left[\int (H-\overline{H})^2\right]^{\frac12}
\end{align}
So
\begin{equation}\label{ine-pp-2}
\int_{\Sigma} |\nabla \phi|^2 \leq  \frac{1}{\eta_1}\int (H-\overline{H})^2.
\end{equation}
Substitute (\ref{ine-pp-2}) and (\ref{eq-p-2}) into (\ref{ine-pp-1}). We obtain 
\begin{align}
&\int_{\Sigma}|H-\overline{H}|^2\nonumber\\
&\leq\frac{n}{n-1}\left(\int_{\Sigma}|\text{\r A}|^2\right)^{\frac12}\left[\frac{n-1}{n}\int_{\Sigma}(H-\overline{H})^2+\left(\frac{(n-1)K}{\eta_1}\right)\int (H-\overline{H})^2\right]^{\frac12}\nonumber\\
&=\sqrt{ \frac{n}{n-1}}\left(\int_{\Sigma}| \text{\r A}|^2\right)^{\frac12}\left[(1+\frac{nK}{\eta_1})\int (H-\overline{H})^2\right]^{\frac12}
\end{align}
So we obtain inequality (\ref{ine-i-04}):
\begin{align}
\int_{\Sigma}|H-\overline{H}|^2\leq \frac{n}{n-1}(1+\frac{nK}{\eta_1})\int_{\Sigma}| \text{\r A}|^2.\nonumber
\end{align}
By (\ref{ine-i-04}) and the identity: $|A-\frac{H}{n}g|^2=|A-\frac{\overline{H}}{n}g|^2-\frac1n(H-\overline{H})^2,$  we obtain (\ref{ine-i-4}):
\begin{equation}\int_{\Sigma}|A-\frac{\overline{H}}{n} g|^2\leq \frac{n}{n-1}(1+\frac{K}{\eta_1})\int_{\Sigma}|A-\frac{H}{n}g|^2.\nonumber
\end{equation}

\bigskip

 Now we prove the conclusion about the equalities in (\ref{ine-i-4}) and  (\ref{ine-i-04}).

Assume  the ambient manifold $M$ is  the Euclidean space $\mathbb{R}^{n+1}$,   the hyperbolic space $\mathbb{H}^{n+1}(-1)$, or  the closed hemisphere $\mathbb{S}_{+}^{n+1}(1)$. Firstly, suppose $\Sigma$ is totally umbilical,  that is, $A=\frac{H}{n}g$. Then the right side   of  (\ref{ine-i-4})  vanishes and the equality holds.  
Secondly,  suppose the equality  in (\ref{ine-i-4}) holds. Thus the equality  in (\ref{ine-i-04}) also holds. We discuss two cases: constant $K=0$ and $K>0$   separately.

{ Case 1.}  When  $K=0$, by the proof of (\ref{ine-i-4}), it is holds that, on $\Sigma$,
\begin{itemize}
\item [(i)] $\text{Ric}(\nabla\phi,\nabla\phi)=0$ and
\item [(ii)] $\text{\r A}$ and $\nabla^2\phi-\frac1n(\Delta\phi)g$ are linearly dependent. 
\end{itemize}

If $\text{\r A}\equiv 0$ on $\Sigma$, $\Sigma$ is totally  umbilical. If $\nabla^2\phi-\frac1n(\Delta\phi)g\equiv 0$ on $\Sigma$,  by (\ref{ine-p-1}), $H=\overline{H}$ on $\Sigma$.  Since   the equality in (\ref{ine-i-04}) holds,   $A=\frac{H}n{g}$ on $\Sigma$.   So  $\Sigma$ is  totally umbilical. If both $\text{\r A}$ and $\nabla^2\phi-\frac1n(\Delta\phi)g$ are not identically zero on $\Sigma$, then  by (ii),  there exists a nonzero constant $\mu\neq 0$ such that, on $\Sigma$,
\begin{align}\label{eq-iii} \mu \text{\r A}=\nabla^2\phi-\frac1n(\Delta\phi)g.
\end{align}

Fix a point  $o\in M$.  Since $\Sigma$ is closed, there is  a point $p\in \Sigma$ such that $p$ realizes the maximum $d$ of the extrinsic distances between a point in $\Sigma$ and $o$ in the metric $\widetilde{g}$ of $M$. Let $B_d(o)$ denote  the closed geodesic ball of $M$ with the radius $d$ centered at $o$.  Then $\Sigma$ is  contained in $B_d(o)$. Since $M^{n+1}(c)=\mathbb{R}^{n+1}, \mathbb{H}^{n+1}, \mathbb{S}_+^{n+1}$,    the distance sphere $S_d(o)=\partial B_d(o)$ is a smooth closed hypersurface and  tangent to $\Sigma$ at $p$. By the Gauss equation, we have 
 \begin{align}\label{eq-gauss}\textrm{Ric}=(n-1)cg+HA-A^2.
 \end{align}
We may choose  $e_i,  1\leq i\leq  n,$ such that they are   the orthonormal  eigenvectors of $A$ at $p$ with $Ae_i=\lambda_ie_i, 1\leq i\leq  n.$  Here we still use $A$ to denote the shape operator of $\Sigma$:  $T\Sigma \rightarrow T\Sigma$,  defined by $\left<A(X),Y\right>= A(X,Y)$. Then
$$\textrm{Ric}(e_i)=c(n-1)e_i+(\displaystyle\sum_{j\neq i}\lambda_j)\lambda_ie_i=[c(n-1)+(\displaystyle\sum_{j\neq i}\lambda_j)\lambda_i]e_i=\tau_ie_i,$$
where  $\tau_i=c(n-1)+(\displaystyle\sum_{j\neq i}\lambda_j)\lambda_i, 1\leq i\leq n.$
This says at $p$,  $e_i$ are the orthonormal  eigenvectors of Ricci tensor,  corresponding to the eigenvalues  $\tau_i$.
  Now we claim that the Ricci tensor $\text{Ric}$ of $\Sigma$ is positive definite at $p$.

Recall  that  the principle curvatures of a hypersurface are the eigenvalues  of its shape operator $A$. It is  known  that  under the above notations, at $p$, the principal curvatures $\lambda_i$   of $\Sigma$  are no less than the principal curvatures $\eta$ of $S_d(o)$ when $M$ is as above. For completeness of proof, we give the proof of this conclusion here.

Fix $i$, $1\leq i\leq n.$ Let $\gamma(s): (-\epsilon, \epsilon)\rightarrow \Sigma$ be a smooth curve with the arc-length parameter, satisfying $\gamma(0)=p, \gamma'(0)=e_i$. Take $h(s)=r(\gamma(s))$, where $r$ denotes the extrinsic distance function from the point $o$. We have $h(0)=r(p)=\displaystyle\max_{s\in (-\epsilon,\epsilon)}h(s)$. Hence $h'(0)=0, h''(0)\leq 0$. Note
$$h'(s)=\left<\widetilde{\nabla}r, \gamma'\right>(s),$$
\begin{align}h''(0)&=\left<\widetilde{\nabla}_{\gamma'(0)}\widetilde{\nabla}r, \gamma'(0)\right>+\left<\widetilde{\nabla}r, \widetilde{\nabla}_{\gamma'(0)}\gamma'\right>(0)\nonumber\\
&=\text{Hess} r(\gamma'(0),\gamma'(0))+\left<\nu,\overline{\nabla}_{\gamma'(0)}\gamma'\right>(0)\nonumber\\
&=\text{Hess} r(e_i,e_i)-A(e_i,e_i)=\eta-\lambda_i\leq 0.\nonumber
\end{align}
Thus we have proved the conclusion  mentioned above.
 
 If  $M = \mathbb{R}^{n+1}$,  then $\eta=\frac{1}{d}$.  So $\lambda_i\geq \eta>0$ and hence  $\tau_i\geq (n-1)\frac{1}{d^2}>0$ at $p$.
 
If $M = \mathbb{H}^{n+1}(-1)$,  it is known  the principle curvature $\eta=\coth d>1$. Then  at $p$, $\lambda_i\geq \eta>1$ and  
 $$\tau_i=-(n-1)+(\displaystyle\sum_{j\neq i}\lambda_j)\lambda_i>0.$$
 
 If $M = \mathbb{S}_{+}^{n+1}(1)$,   the principle curvature $\eta=\cot d\geq 0$. Then  at $p$, $\lambda_i\geq \eta\geq 0$ and
  $$\tau_i=(n-1)+(\displaystyle\sum_{j\neq i}\lambda_j)\lambda_i>0.$$

This means that  the Ricci tensor  of $\Sigma$ is positive definite at $p$. We have proved the claim. 

Then there exists a neighborhood $N_{\epsilon}( p )$ of $p$ in $\Sigma$ such that the Ricci tensor of $\Sigma$ is positive-definite on $N_{\epsilon}( p )$. 
Note  (i)  $\text{Ric}(\nabla\phi,\nabla\phi)=0.$  It must hold that $\nabla \phi=0$ on $N_{\epsilon}( p )$. By this and the definition of $\phi$,   we have $\phi\equiv C$ and $H=\overline{H}$ on $N_{\epsilon}( p )$. By (\ref{eq-iii}),  $A=\frac{H}{n}g=\frac{\overline{H}}{n}g$  on $N_{\epsilon}( p )$.  By  continuity, on the closed neighborhood  $\overline{N}_{\epsilon}( p )$,  $A=\frac{\overline{H}}{n}g$. 
Hence by (\ref{eq-gauss}), $\textrm{Ric}=(n-1)cg+\frac{\overline{H}^2}{n}g-\frac{\overline{H}^2}{n^2}g^2$ on $\overline{N}_{\epsilon}( p )$, which says that on $\overline{N}_{\epsilon}( p )$,  the Ricci tensor is a constant tensor and thus is positive definite.  

By the same argument as above,  we may have   a neighborhood of every point on the boundary $\partial N_{\epsilon}( p )$ of $N_{\epsilon}( p )$ such that $A=\frac{\overline{H}}{n}g$. Let $D$ be the connected subset of $\Sigma$ such that $p\in D$ and $A=\frac{\overline{H}}{n}g$. From the above argument,  $D$ is both open and closed. Since $\Sigma$ is connected,  $D=\Sigma$ and thus $\Sigma$ is totally umbilical.
 
{Case 2.} When  $K>0$,  we have that, on $\Sigma$,
\begin{itemize}
\item [(I)] $\left(\text{Ric}+(n-1)Kg\right)(\nabla \phi,\nabla \phi)=0,$
\item [(II)] $\text{\r A}$ and $\nabla^2\phi-\frac1n(\Delta \phi)g$ are linearly dependent,
\item[(III)] $\phi$ and $H-\overline{H}$ are linearly dependent, and
\item [(IV)]  $\int_{\Sigma}(|\nabla \phi|^2-\eta_1\phi^2)=0$.
\end{itemize}
Like the proof in Case 1, $M$ is obviously totally umbilical if  $\text{\r A}\equiv 0$ ou $\nabla^2\phi-\frac1n(\Delta\phi)g\equiv 0$.  We now consider that  both $\text{\r A}$ and $\nabla^2\phi-\frac1n(\Delta\phi)g$ are not identically zero on $\Sigma$, then  by (II),  there exists a nonzero constant $\mu\neq 0$ such that

\begin{align}\label{eq-k-1} \mu \text{\r A}=\nabla^2\phi-\frac1n(\Delta\phi)g.
\end{align}
Take points $o$ and $p$ as in  Case 1. Observe that  the positivity of Ricci curvature on  the neighborhood $N_{\epsilon}( p )$ of $p$ and $K>0$ still  implies $\nabla \phi=0$ on  $N_{\epsilon}( p )$.  Similar to case $1$, we may prove that $\Sigma$ is totally umbilical.

We complete the proof of the theorem.
\qed

 Theorem \ref{thm-3} is obtained by taking  $K=0$ in Theorem \ref{thm-4}.

\bigskip
\bigskip
\noindent  Xu Cheng\\Insitituto de Matem\'atica\\Universidade
Federal Fluminense - UFF\\Centro, Niter\'{o}i, RJ 24020-140 Brazil
\\e-mail:xcheng@impa.br

\bigskip
\bigskip
\noindent Detang Zhou\\Insitituto de Matem\'atica\\Universidade
Federal Fluminense - UFF\\Centro, Niter\'{o}i, RJ 24020-140 Brazil
\\e-mail: zhou@impa.br
\end{document}